\documentclass[]{interact}

\usepackage[utf8]{inputenc}
\usepackage[T1]{fontenc}
\usepackage[english]{babel}
\usepackage{lmodern}
\usepackage{amsmath}
\usepackage{amssymb}
\usepackage{mathrsfs}
\usepackage{amsthm}
\usepackage{color}
\usepackage{soul}
\usepackage{diagbox}
\usepackage{ulem}
\usepackage{verbatim}
\usepackage{moreverb}
\usepackage{url}
\usepackage{graphicx}
\usepackage{listings}
\usepackage{dsfont}
\usepackage[breaklinks]{hyperref}
\usepackage{layout}
\hypersetup{urlcolor=blue,linkcolor=blue,citecolor=blue,colorlinks=true}

\newtheorem{theo}{Theorem}[section]
\newtheorem{defi}{Definition}[subsection]

\newtheorem{Rem}{Remark}[section]
\newtheorem{Cor}{Corollary}[section]
        
\begin{document}


\title{Asymptotic distribution of least square estimators for linear models with dependent errors}

\author{
\name{Emmanuel Caron\textsuperscript{a}\thanks{Corresponding author. Email: emmanuel.caron@ec-nantes.fr}}
\affil{\textsuperscript{a} Ecole Centrale Nantes, Laboratoire de Mathématiques Jean Leray UMR 6629, 1 Rue de la Noë, 44300 Nantes, France}
}

\maketitle

%
%
%

\pagestyle{plain}

\begin{abstract}
In this paper, we consider the usual linear regression model in the case where the error process is assumed strictly stationary. We use a result from Hannan \cite{hannan73clt}, who proved a Central Limit Theorem for the usual least square estimator under general conditions on the design and on the error process. 
Whatever the design satisfying Hannan's conditions, we define an estimator of the covariance matrix and we prove its consistency under very mild conditions. As an application, we show how to modify the usual tests on the linear model in this dependent context, in such a way that the type-$I$ error rate remains asymptotically correct, and we illustrate the performance of this procedure through different sets of simulations.
\end{abstract}


\section{Introduction}

The linear regression model is used in many domains of applied mathematics, and the asymptotic behavior of the least square estimators is well known when the errors are i.i.d. (independent and identically distributed) random variables.
Many authors have deepened the research on this subject, we can cite for example~\cite{bassett1978asymptotic},~\cite{babu1989strong},~\cite{bai1992m} and~\cite{he1996general} among others. 
However, many science and engineering data exhibit significant temporal dependence so that the assumption of independence is violated (see for instance~\cite{brockwell2013time}). It is observed in astrophysics, geophysics, biostatistics, climatology, among others. Consequently all statistical procedures based on this assumption are not efficient and this can be very problematic for the applications.

In this paper, we propose to study the usual linear regression model in the very general framework of Hannan~\cite{hannan73clt}. Let us consider the equation of the model: 
\[Y = X\beta + \epsilon.\]
The process $(\epsilon_{i})_{i \in \mathbb{Z}}$ is assumed to be strictly stationary. The $n \times p$ matrix $X$ is the design and can be random or deterministic. In our framework, we consider the inter-dependence of the variables of the design.
As in Hannan, we assume that the design matrix $X$ is independent of the error process.
Such a model can be used for time series regression, but also in a more general context when the residuals seem to derive from a stationary correlated process. 

Our work is based on the paper by Hannan \cite{hannan73clt}, who proved a Central Limit Theorem for the usual least square estimator under general conditions on the design and on the error process. Let us quote that most of short-range dependent processes satisfies Hannan’s conditions on the error process, for instance the class of linear processes with summable coefficients and square integrable innovations, a large class of functions of linear processes, many processes under various mixing conditions and the $2$-strong stable processes introduced by Wu~\cite{wu2005nonlinear}.
We refer to our previous paper~\cite{caron1}, which presents many classes of short-range dependent processes satisfying Hannan’s condition.

The linear regression model with dependent errors has also been studied under more restrictive conditions. For instance, Pagan and Nicholls~\cite{pagan1976exact} consider the case where the errors follow a $MA(q)$ process, and Chib and Greenberg~\cite{chib1994bayes} the case where the errors are an $ARMA(p,q)$ process.
A more general framework is used by Wu~\cite{wu2007m} for a class of short-range dependent processes. These results are based on the asymptotic theory of stationary processes developed by Wu in~\cite{wu2005nonlinear}. However the class of processes satisfying the so called $\mathbb{L}^{2}$ "physical dependence measure" is included in the class of processes satisfying Hannan's condition~\eqref{C1}. 
In the present paper, we consider the very general framework of Hannan in order to obtain the most robust results. 

In this paper, we present an estimator of the asymptotic covariance matrix of the normalized least square estimators of the parameters. This estimator is derived from the estimator of the spectral density of the error process introduced in~\cite{caron1}.
Once the asymptotic covariance matrix is consistently estimated, it is then possible to obtain confidence regions and test procedures for the unknown parameter $\beta$. 
In particular, we shall use our general results to modify the usual Student and Fisher tests in cases where $(\epsilon_{i})_{i \in \mathbb{Z}}$ and the design verify the conditions of Hannan, in order to have always a type-$I$ error rate asymptotically correct (approximately equals to $5$\%). 

The paper is organized as follows. In Section $2$, we recall Hannan’s Central Limit Theorem for the least square estimator.
In Section $3$, we focus on the estimation of the covariance matrix under Hannan's conditions. 
Finally, Section $4$ is devoted to the correction of the usual Student and Fisher tests in our dependent context, and some simulations with different models are realized.

\section{Hannan's theorem}

\subsection{Notations and definitions}

Let us recall the equation of the linear regression model: 
\begin{equation}
Y = X\beta + \epsilon,
\label{mod_lin}
\end{equation}
where $X$ is a design matrix and $\epsilon$ is an error process defined on a probability space ($\Omega, \mathcal{F}, \mathbb{P}$). Let us notice that the error process $\epsilon$ is independent of the design $X$. Let $X_{.,j}$ be the column $j$ of the matrix $X$, and $x_{i,j}$ the real number at the row $i$ and the column $j$, where $j$ is in $\{1, \ldots, p\}$ and $i$ in $\{1, \ldots, n\}$. The random vectors $Y$ and $\epsilon$ belong to $\mathbb{R}^{n}$ and $\beta$ is a $p \times 1$ vector of unknown parameters.

Let $\left \| . \right \|_{2}$ be the usual euclidean norm on $\mathbb{R}^{n}$, and $\left \| . \right \|_{\mathbb{L}^{p}}$ be the $\mathbb{L}^{p}$-norm on $\Omega$, defined for all random variable $Z$ by: $\left \| Z \right \|_{\mathbb{L}^{p}} = \left[ \mathbb{E} \left( Z^{p} \right) \right]^{\frac{1}{p}}$ . We say that $Z$ is in $\mathbb{L}^{p}(\Omega)$ if $\left[ \mathbb{E} \left( Z^{p} \right) \right]^{\frac{1}{p}} < \infty$.

The error process $(\epsilon_{i})_{i \in \mathbb{Z}}$ is assumed to be strictly stationary with zero mean. Moreover, for all $i$ in $\mathbb{Z}$, $\epsilon_{i}$ is supposed to be in $\mathbb{L}^{2}(\Omega)$. More precisely, the error process satisfies, for all $i$ in $\mathbb{Z}$: 
\[\epsilon_{i} = \epsilon_{0} \circ \mathbb{T}^{i},\]
where $\mathbb{T}: \Omega \rightarrow \Omega$ is a bijective bimeasurable transformation  preserving the probability measure $\mathbb{P}$. Note that any strictly stationary process can be represented in this way.

Let ($\mathcal{F}_{i}$)$_{i \in \mathbb{Z}}$ be a non-decreasing filtration built as follows, for all $i$: 
\[\mathcal{F}_{i} = \mathbb{T}^{-i}(\mathcal{F}_{0}),\]
where $\mathcal{F}_{0}$ is a sub-$\sigma$-algebra of $\mathcal{F}$ such that $\mathcal{F}_{0} \subseteq \mathbb{T}^{-1}(\mathcal{F}_{0})$. For instance, one can choose the past $\sigma$-algebra before time $0$: $\mathcal{F}_{0} = \sigma(\epsilon_{k}, k \leq 0)$, and then $\mathcal{F}_{i} = \sigma(\epsilon_{k}, k \leq i)$. In that case, $\epsilon_{0}$ is $\mathcal{F}_{0}$-measurable.

As in Hannan, we shall always suppose that $\mathcal{F}_{-\infty} = \underset{i \in \mathbb{Z}}{\bigcap} \mathcal{F}_{i}$ is trivial. Moreover $\epsilon_{0}$ is assumed $\mathcal{F}_{\infty}$-measurable. These implie that the $\epsilon_{i}$'s are all regular random variables in the following sense:

\begin{defi}[Regular random variable]
Let $Z$ be a random variable in $L^{1}(\Omega)$. We say that $Z$ is regular with respect to the filtration $(\mathcal{F}_{i})_{i \in \mathbb{Z}}$ if $\mathbb{E}(Z | \mathcal{F}_{-\infty}) = \mathbb{E}(Z)$ almost surely and if $Z$ is $\mathcal{F}_{\infty}$-measurable. 
\end{defi}
Hence there exists a spectral density $f$ for the error process, defined on $[-\pi, \pi]$. The autocovariance function $\gamma$ of the process $\epsilon$ then satisfies: 
\[\gamma(k) = \mathrm{Cov} (\epsilon_{m}, \epsilon_{m+k}) = \mathbb{E}(\epsilon_{m}\epsilon_{m+k}) = \int_{-\pi}^{\pi} e^{i k \lambda} f(\lambda) d\lambda.\]
Furthermore we denote by $\Gamma_{n}$ the covariance matrix of the error process:
\begin{equation}
\Gamma_{n} = \left[ \gamma(j-l) \right]_{1 \leq j,l \leq n}.
\label{Grand_Gamma}
\end{equation}

\subsection{Hannan's Central Limit Theorem}

Let $\hat{\beta}$ be the usual least square estimator for the unknown vector $\beta$.
Given the design $X$, Hannan \cite{hannan73clt} has shown a Central Limit Theorem for $\hat{\beta}$ when the error process is stationary. In this section, the conditions for applying this theorem are recalled.

Let $(P_{j})_{j \in \mathbb{Z}}$ be a family of projection operators, defined for all $j$ in $\mathbb{Z}$ and for any $Z$ in $\mathbb{L}^{2}(\Omega)$ by:
\[P_{j}(Z) = \mathbb{E}(Z | \mathcal{F}_{j}) - \mathbb{E}(Z | \mathcal{F}_{j-1}).\]
We shall always assume that Hannan's condition on the error process is satisfied:
\begin{equation}
\sum_{i \in \mathbb{Z}} \left \| P_{0}(\epsilon_{i}) \right \|_{\mathbb{L}^{2}} < +\infty.
\tag{C1}
\label{C1}
\end{equation}
Note that this condition implies that:
\begin{equation}
\sum_{k \in \mathbb{Z}} \left| \gamma(k) \right| < \infty,
\label{0bis}
\end{equation}
(see for instance \cite{dmv2007weak}).

Hannan's condition provides a very general framework for stationary processes. The hypothesis~\eqref{C1} is a sharp condition to have a Central Limit Theorem for the partial sum sequence (see the paper of Dedecker, Merlevède and Voln\'y \cite{dmv2007weak} for more details). Notice that the condition~\eqref{0bis} implies that the error process is short-range dependent.
However, Hannan's condition is satisfied for most short-range dependent stationary processes. The reader can see the paper~\cite{caron1} where some examples checking Hannan's condition are developed.\\

Let us now recall Hannan’s assumptions on the design. Let us introduce:
\[d_{j}(n) = \left \| X_{.,j} \right \|_{2} = \sqrt{\sum_{i=1}^{n} x_{i, j}^{2}},\]
and let $D(n)$ be the diagonal matrix with diagonal term $d_{j}(n)$ for $j$ in $\{1, \ldots, p\}$.

Following Hannan, we also require that the columns of the design $X$ satisfy, almost surely, the following conditions:

\begin{equation}
\forall j \in \{1, \ldots,p\}, \qquad \lim_{n \rightarrow \infty} d_{j}(n) = \infty,
\tag{C2}
\label{C2}
\end{equation}
and:
\begin{equation}
\forall j \in \{1, \ldots, p\}, \qquad \lim_{n \rightarrow \infty} \frac{\sup_{1 \leq i \leq n} \left | x_{i,j} \right |}{d_{j}(n)} = 0.
\tag{C3}
\label{C3}
\end{equation}
Moreover, we assume that the following limits exist: 

\begin{equation}
\forall j, l \in \{1, \ldots, p\}, \ k \in \{0, \ldots, n-1\}, \qquad \rho_{j,l}(k) = \lim_{n \rightarrow \infty} \sum_{m=1}^{n-k} \frac{x_{m, j} x_{m+k,l}}{d_{j}(n)d_{l}(n)}.
\tag{C4}
\label{C4}
\end{equation}

Note that Conditions~\eqref{C2} and~\eqref{C3} correspond to the usual Lindeberg’s conditions for linear statistics in the i.i.d. case. In the dependent case, we also need Condition~\eqref{C4}.

The $p \times p$ matrix formed by the coefficients $\rho_{j,l}(k)$ is called $R(k)$:
\begin{equation}
R(k) = [\rho_{j,l}(k)] = \int_{-\pi}^{\pi} e^{i k \lambda} F_{X}(d\lambda), \quad a.s.
\label{4}
\end{equation}
where $F_{X}$ is the spectral measure associated with the matrix $R(k)$. The matrix $R(0)$ is supposed to be positive definite: 
\begin{equation}
R(0) > 0, \quad a.s.
\tag{C5}
\label{C5}
\end{equation}
Let then $F$ and $G$ be the matrices:
\[F = \frac{1}{2\pi} \int_{-\pi}^{\pi} F_{X}(d\lambda), \quad a.s.\]
\[G = \frac{1}{2\pi} \int_{-\pi}^{\pi} F_{X}(d\lambda) \otimes f(\lambda), \quad a.s.\]

The Central Limit Theorem for the regression parameter, due to Hannan \cite{hannan73clt}, can be stated as follows:

\begin{theo}
Let $(\epsilon_{i})_{i \in \mathbb{Z}}$ be a stationary process with zero mean. Assume that $\mathcal{F}_{-\infty}$ is trivial, $\epsilon_{0}$ is $\mathcal{F}_{\infty}$-measurable, and that the sequence $(\epsilon_{i})_{i \in \mathbb{Z}}$ satisfies Hannan's condition~\eqref{C1}. Assume that the design $X$ satisfies, almost surely, the conditions~\eqref{C2}, \eqref{C3}, \eqref{C4} and~\eqref{C5}. 
Then, for all bounded continuous function $f$: 
\begin{equation}
\mathbb{E} \left( f \left( D(n)(\hat{\beta} - \beta) \right) \Big| X \right) \xrightarrow[n \rightarrow \infty]{} \mathbb{E} \left( f(Z) \Big| X \right), \quad a.s.
\label{7}
\end{equation}
where the distribution of $Z$ given $X$ is a gaussian distribution, with mean zero and covariance matrix equal to $F^{-1}GF^{-1}$.
Furthermore, there is the convergence of the second order moment: \footnote{The transpose of a matrix $X$ is denoted by $X^{t}$.}
\begin{equation}
\mathbb{E} \left( D(n) (\hat{\beta} -\beta) (\hat{\beta} -\beta)^{t} D(n)^{t} \Big| X \right) \xrightarrow[n \rightarrow \infty]{} F^{-1}GF^{-1}, \quad a.s.
\label{9}
\end{equation}
\label{Hannan_th}
\end{theo}

\begin{Rem}
Let us notice that, by the dominated convergence theorem, the property~\eqref{7} implies that for any bounded continuous function $f$, 
\[\mathbb{E} \left( f \left( D(n)(\hat{\beta} - \beta) \right) \right) \xrightarrow[n \rightarrow \infty]{} \mathbb{E} \left( f(Z) \right).\]
\end{Rem}

\begin{Rem}
In this remark, for the sake of clarity, we give a direct proof of~\eqref{9}. We shall see that, in fact,~\eqref{9} holds under~\eqref{0bis} and~\eqref{C4} - \eqref{C5} (Hannan's condition~\eqref{C1}, which implies~\eqref{0bis}, is needed for~\eqref{7} only). Moreover, this proof will serve as a preliminary to the proof of Theorem~\ref{main_th}.
We start from the exact expression of the second order moment:
\begin{multline*}
\mathbb{E} \left( D(n) (\hat{\beta} -\beta) (\hat{\beta} -\beta)^{t} D(n)^{t} \Big| X \right) \\
= D(n) (X^{t}X)^{-1} X^{t} \Gamma_{n} X (X^{t}X)^{-1} D(n) \\
= \hat{R}(0)^{-1} \left( D(n)^{-1} X^{t} \Gamma_{n} X D(n)^{-1} \right) \hat{R}(0)^{-1},
\end{multline*}
with $\hat{R}(0) = D(n)^{-1} X^{t} X D(n)^{-1}$. 
The $n \times n$ covariance matrix $\Gamma_{n}$ is a symmetric Toeplitz matrix and is equal to:
\[\Gamma_{n} = \sum_{k=-n+1}^{n-1} \gamma(k) J_{n}^{(k)},\]
where $J^{(k)}$ is the matrix with some $1$ on the $k$-th diagonal and $0$ elsewhere.

Hence, we deduce that:
\[\mathbb{E} \left( D(n) (\hat{\beta} -\beta) (\hat{\beta} -\beta)^{t} D(n)^{t} \Big| X \right) = \hat{R}(0)^{-1} \left( \sum_{k=-n+1}^{n-1} \gamma(k) B_{k,n} \right) \hat{R}(0)^{-1},\]
with:
\[B_{k,n} = D(n)^{-1} X^{t} J_{n}^{(k)} X D(n)^{-1}.\]
For all $k$ in $\{-n+1, \dots, n-1\}$, the matrices $B_{k,n}$ are equal to:
\begin{equation}
B_{k,n} = [\hat{\rho}_{j,l}(k)] \quad \text{if} \ k \geq 0, \qquad B_{k,n} = [\hat{\rho}_{j,l}(-k)] \quad \text{if} \ k \leq -1,
\label{Bkn}
\end{equation}
where $\hat{\rho}_{j,l}(k) = \sum_{m=1}^{n-k} \frac{x_{m, j} x_{m+k,l}}{d_{j}(n)d_{l}(n)}$. Under~\eqref{C4}, $\hat{\rho}_{j,l}(k)$ converges almost surely to $\rho_{j,l}(k)$.

By the dominated convergence theorem, since every term of $B_{k,n}$ is in $[-1,1]$, we deduce that:
\[\sum_{k=-n+1}^{n-1} \gamma(k) B_{k,n} \xrightarrow[n \rightarrow \infty]{} \sum_{k=-\infty}^{\infty} \gamma(k) B_{k},\]
where $B_{k} = [\rho_{j,l}(k)]$ if $k \geq 0$ and $B_{k} = [\rho_{j,l}(-k)]$ if $k \leq -1$.

Since moreover $\hat{R}(0)$ converges almost surely to $R(0)$ (which is positively definite, see ~\eqref{C5}) as $n$ tends to infinity, we conclude that:
\[\mathbb{E} \left( D(n) (\hat{\beta} -\beta) (\hat{\beta} -\beta)^{t} D(n)^{t} \Big| X \right) \xrightarrow[n \rightarrow \infty]{} R(0)^{-1} \left( \sum_{k=-\infty}^{\infty} \gamma(k) B_{k} \right) R(0)^{-1}.\]

Note that $R(0) = \int_{-\pi}^{\pi} F_{X} (d \lambda) = 2 \pi F$ and $\sum_{k=-\infty}^{\infty} \gamma(k) B_{k} = 4 \pi^{2} G$, which is consistent with~\eqref{9}.
\label{Rem2.2}
\end{Rem}

\section{Estimation of the covariance matrix}

To obtain confidence regions or test procedures from Theorem~\ref{Hannan_th}, one needs to estimate the limiting covariance matrix $F^{-1} G F^{-1}$. In this section, we propose an estimator of this covariance matrix, and we show its consistency under Hannan's conditions. 

Let us first consider a preliminary random matrix defined as follows:
\begin{equation}
\widehat{\Gamma}_{n, h_{n}} = \left[ K \left( \frac{j-l}{h_{n}} \right) \hat{\gamma}_{j-l} \right]_{1 \leq j,l \leq n},
\label{Gamma_tapered}
\end{equation}
with:
\[\hat{\gamma}_{k} = \frac{1}{n} \sum_{j=1}^{n-|k|} \epsilon_{j} \epsilon_{j+|k|}, \quad 0 \leq | k | \leq n-1.\]
The function $K$ is a kernel such that:
\begin{itemize}
\item[-] $K$ is nonnegative, symmetric, and $K(0) = 1$,
\item[-] $K$ has compact support,
\item[-] The Fourier transform of $K$ is integrable.
\end{itemize}
The sequence of positive reals $h_{n}$ is such that $h_{n}$ tends to infinity and $\frac{h_{n}}{n}$ tends to $0$ when $n$ tends to infinity.

In our context, the errors $(\epsilon_{i})_{1 \leq i \leq n}$ are not observed. Only the residuals are available:
\[\hat{\epsilon}_{j} = Y_{j} - (x_{j})^{t} \hat{\beta},\]
because only the data $Y$ and the design $X$ are observed. Consequently, we consider the following estimator of $\Gamma_{n}$:
\begin{equation}
\widehat{\Gamma}_{n, h_{n}}^{\ast} = \left[ K \left( \frac{j-l}{h_{n}} \right) \hat{\gamma}_{j-l}^{\ast} \right]_{1 \leq j,l \leq n},
\label{Gamma_tapered_star}
\end{equation}
with:
\[\hat{\gamma}_{k}^{\ast} = \frac{1}{n} \sum_{j=1}^{n-|k|} \hat{\epsilon}_{j} \hat{\epsilon}_{j+|k|}, \quad 0 \leq | k | \leq n-1.\]

This estimator is a truncated version of the full matrix $\widehat{\Gamma}_{n}^{\ast} = \left[ \hat{\gamma}_{j-l}^{\ast} \right]_{1 \leq j,l \leq n}$, preserving the diagonal and some sub-diagonals. Following~\cite{bickel2008regularized}, $\widehat{\Gamma}_{n, h_{n}}^{\ast}$ is called the tapered covariance matrix estimator. The motivation for tapering comes from the fact that, for a large $k$, either $\gamma(k)$ is close to zero or $\hat{\gamma}_{k}^{\ast}$ is an unreliable estimate of $\gamma(k)$. Thus, prudent use of tapering may bring considerable computational economy in the former case, and statistical efficiency in the simulations, by keeping small or unreliable $\hat{\gamma}_{k}^{\ast}$ out of the calculations.

To estimate the asymptotic covariance matrix $F^{-1} G F^{-1}$, we use the estimator:
\begin{equation}
C_{n} = D(n) (X^{t}X)^{-1} X^{t} \widehat{\Gamma}_{n, h_{n}}^{\ast} X (X^{t}X)^{-1} D(n).
\label{Cov_estimate}
\end{equation}
Let us denote by $C$ the matrix $F^{-1} G F^{-1}$ and the coefficients of the matrices $C_{n}$ and $C$ are respectively denoted by $c_{n,(j,l)}$ and $c_{j,l}$, for all $j, l$ in $\{1, \ldots, p\}$. Our first result is the following:

\begin{theo}
Let $h_{n}$ be a sequence of positive reals such that $h_{n} \rightarrow \infty$ as $n$ tends to infinity, and:
\begin{equation}
h_{n} \mathbb{E} \left( \left| \epsilon_{0} \right|^{2} \left( 1 \wedge \frac{h_{n}}{n} \left| \epsilon_{0} \right|^{2} \right) \right) \xrightarrow[n \rightarrow \infty]{} 0.
\label{c_n_cond}
\end{equation}
Then, under the assumptions of Theorem~\ref{Hannan_th}, the estimator $C_{n}$ is consistent, that is for all $j, l$ in $\{1, \ldots, p\}$:
\begin{equation}
\mathbb{E} \left( \left| c_{n,(j,l)} - c_{j,l} \right| \Big| X \right) \xrightarrow[n \rightarrow \infty]{} 0, \quad \text{a.s.}
\label{consistence_eq}
\end{equation}
\label{main_th}
\end{theo}

\begin{Rem}
If $\epsilon_{0}$ is square integrable, then there exists $h_{n} \rightarrow \infty$ such that~\eqref{c_n_cond} holds.

Furthermore if $\mathbb{E} \left( \left | \epsilon_{0} \right |^{\delta+2} \right) < \infty$ with $\delta$ in $]0,2]$, then:
\[h_{n} \mathbb{E} \left( \left| \epsilon_{0} \right|^{2} \left( 1 \wedge \frac{h_{n}}{n} \left| \epsilon_{0} \right|^{2} \right) \right)  \leq h_{n} \mathbb{E} \left( \left| \epsilon_{0} \right|^{2} \frac{h_{n}^{\delta/2}}{n^{\delta/2}} |\epsilon_{0}|^{\delta} \right) \leq \frac{h_{n}^{1+\delta/2}}{n^{\delta/2}} \mathbb{E} \left( \left | \epsilon_{0} \right |^{\delta+2} \right).\] 
Thus, if $h_{n}$ satisfies $\frac{h_{n}^{1+\delta/2}}{n^{\delta/2}} \xrightarrow[n \rightarrow \infty]{} 0$, then~\eqref{c_n_cond} holds.
In particular, if the random variable $\epsilon_{0}$ has a fourth order moment, then the condition on $h_{n}$ is $\frac{h_{n}^{2}}{n} \xrightarrow[n \rightarrow \infty]{} 0$.
\end{Rem}

From this theorem, we get the non-conditional convergence in probability:
\begin{Cor}
Let $h_n$ be a sequence satisfying~\eqref{c_n_cond}. Then the estimator $C_{n}$ converges in probability to $C$ as $n$ tends to infinity.
\label{cor_main_th}
\end{Cor}

\begin{Rem}
Since $F^{-1} G F^{-1}$ is assumed to be positive definite, our estimator $C_{n}$ is also asymptotically positive definite. But it has no reason to be positive definite for any kernel and for any $n$. To overcome this problem, one can consider the estimator $\tilde C_n$ which is built as $C_{n}$ but with a positive definite kernel, like for instance the triangular kernel.

Indeed, following Wu~\cite{xiao2012covariance}, we can define:
\[\widehat{\Gamma}_{n, h_{n}}^{\ast} = \widehat{\Gamma}_{n}^{\ast} \star W_{n},\]
where $\star$ is the Hadamard (or Schur) product, which is formed by element-wise multiplication of matrices, and $W_{n}$ is the kernel's matrix equal to $\left[ K \left( \frac{j-l}{h_{n}} \right) \right]_{1 \leq j,l \leq p}$. Let us notice that the full matrix $\widehat{\Gamma}_{n}^{\ast}$ is positive definite if and only if $\hat{\gamma}_{0}^{\ast} > 0$ (see~\cite{brockwell2013time}). 
Consequently, by the Schur Product Theorem in matrix theory~\cite{horn1990matrix}, since $\hat{\Gamma}_{n}^{\ast}$ and $W_{n}$ are both positive definite, their Schur product $\widehat{\Gamma}_{n, h_{n}}^{\ast}$ is also positive definite.

Let us recall that $C_{n} = \Psi \widehat{\Gamma}_{n, h_{n}}^{\ast} \Psi^{t}$ with $\Psi = D(n) (X^{t} X)^{-1} X^{t}$. Then the estimator $C_{n}$ is positive definite if for all $x \neq 0$, $x^{t} \Psi \widehat{\Gamma}_{n, h_{n}}^{\ast} \Psi^{t} x$ is strictly greater than $0$. 
It is true if $\hat{\gamma}_{0}^{\ast} > 0$ and if the design $X$ is a rank $p$ matrix.
\end{Rem}

Combining Theorem~\ref{Hannan_th} and Theorem~\ref{main_th}, we obtain the following corollary, which is the main result of our paper:
\begin{Cor}
Under the assumptions of Theorem~\ref{Hannan_th} and Theorem~\ref{main_th}, we get:
\begin{equation}
C_{n}^{-\frac{1}{2}} \left( D(n) (\hat{\beta} - \beta) \right) \xrightarrow[n \rightarrow \infty]{\mathcal{L}} \mathcal{N}(0,I_{p}),
\label{Hannan_estim}
\end{equation}
where $I_{p}$ is the $p \times p$ identity matrix.
\label{TLC_estimator}
\end{Cor}

\section{Tests and simulations}

As an application of this main result, we show how to modify the usual tests on the linear regression model.

\subsection{Tests}

Let us recall the assumptions. We consider the linear regression model~\eqref{mod_lin}, and we assume that Hannan's condition~\eqref{C1} as well as the conditions~\eqref{C2} to~\eqref{C5} on the design are satisfied.
We also assume that $\epsilon_{0}$ is $\mathcal{F}_{\infty}$-measurable and that $\mathcal{F}_{-\infty}$ is trivial.
With these conditions, the usual tests can be modified and adapted to the case where the errors are short-range dependent and for any design verifying Hannan's conditions.

As usual, the null hypothesis $H_{0}$ means that the parameter $\beta$ belongs to a vector space with dimension strictly smaller than $p$, and we denote by $H_{1}$ the alternative hypothesis (meaning that $H_{0}$ is not true, but~\eqref{mod_lin} holds).

In order to test $H_{0}: \beta_{j} = 0$ against $H_{1}: \beta_{j} \neq 0$, for $j$ in $\{1, \ldots, p\}$, under the $H_{0}$-hypothesis and according to Corollary~\ref{TLC_estimator} we have: 
\[d_{j}(n) \hat{\beta}_{j} \xrightarrow[n \rightarrow \infty]{} \mathcal{N}(0, c_{j,j}).\]
We introduce the following univariate test statistic:
\begin{equation}
T_{j,n} = \frac{d_{j}(n) \hat{\beta}_{j}}{\sqrt{c_{n,(j,j)}}}.
\label{Student}
\end{equation}
Under the $H_{0}$-hypothesis, the distribution of $T_{j,n}$ converges to a standard normal distribution when $n$ tends to infinity.\\

Now we want test $H_{0}$: $\beta_{j_{1}} = \ldots = \beta_{j_{p_{0}}} = 0$, against $H_{1}$: $\exists j_{z} \in \{j_{1}, \ldots, j_{p_{0}}\}$ such that $\beta_{j_{z}} \neq 0$. By Corollary~\ref{TLC_estimator}, it follows that:
\[
C_{n_{p_{0}}}^{-1/2}
\begin{pmatrix}
    d_{j_{1}}(n) (\hat{\beta}_{j_{1}} - \beta_{j_{1}}) & \\
    \vdots & \\
    d_{j_{p_{0}}}(n) (\hat{\beta}_{j_{p_{0}}} - \beta_{j_{p_{0}}}) & \\
\end{pmatrix}
\xrightarrow[n \rightarrow \infty]{\mathcal{L}} \mathcal{N}(0_{p_{0} \times 1}, I_{p_{0}}),
\]
where $C_{n_{p_{0}}}$ is the covariance matrix $C_{n}$ built with removing the rows and the columns which do not belong to the discrete set $\{j_{1}, \ldots, j_{p_{0}}\}$. The $p_{0} \times p_{0}$ identity matrix is denoted by $I_{p_{0}}$ and $0_{p_{0} \times 1}$ is a $p_{0}$ vector of zeros. 

Then under $H_{0}$-hypothesis, we have:
\[
\begin{pmatrix}
    Z_{1,n} & \\
    \vdots & \\
    Z_{p_{0},n} & \\
\end{pmatrix}
=
C_{n_{p_{0}}}^{-1/2}
\begin{pmatrix}
    d_{j_{1}}(n) \hat{\beta}_{j_{1}} & \\
    \vdots & \\
    d_{j_{p_{0}}}(n) \hat{\beta}_{j_{p_{0}}} & \\
\end{pmatrix}
\xrightarrow[n \rightarrow \infty]{\mathcal{L}} \mathcal{N}(0_{p_{0} \times 1}, I_{p_{0}}),
\]
and we define the following test statistic:
\begin{equation}
\Xi = Z_{1,n}^{2} + \cdots + Z_{p_{0},n}^{2}.
\label{Fisher}
\end{equation}
Under the $H_{0}$-hypothesis, the distribution $\Xi$ converges to a $\chi^{2}$-distribution with parameter $p_{0}$. \newline

For the simulations, we shall use for the estimator $C_{n}$ the kernel $K$ defined by:
\begin{equation}
\left\{
\begin{array}{r c l}
K(x) &=& 1 \phantom{- |x| 1 1} \quad  if\  |x| < 0.8\\
K(x) &=& 5 - 5|x| \phantom{\quad} if\  0.8 \leq |x| \leq 1\\
K(x) &=& 0 \phantom{- |x| 1 1}  \quad if\  |x| > 1.\\
\end{array}
\right.
\label{kernel}
\end{equation}

This kernel verifies the conditions defined at the beginning of Section $3$,
and it is close to the rectangular kernel (whose Fourier transform is not integrable).
Hence, the parameter $h_{n}$ can be understood as the number of covariance terms that are necessary to obtain a good approximation of $\Gamma_{n}$.
To choose its values, we shall use the graph of the empirical autocovariance of the residuals. 

\subsection{Simulations}

We first simulate $(Z_{1}, \ldots, Z_{n})$ according to the $AR(1)$ equation $Z_{k+1} = \frac{1}{2} (Z_{k} + \eta_{k+1})$, where $Z_{1}$ is uniformly distributed over $[0,1]$ and $(\eta_{i})_{i \geq 2}$ is a sequence of i.i.d. random variables with distribution $\mathcal{B}(1/2)$, independent of $Z_{1}$. 
The transition kernel of the chain $(Z_{i})_{i \geq 1}$ is: 
\[K(f)(x) = \frac{1}{2} \left( f \left( \frac{x}{2} \right) + f \left( \frac{x+1}{2} \right) \right),\]
and the uniform distribution on $[0,1]$ is the unique invariant distribution by $K$. Hence, the chain $(Z_{i})_{i \geq 1}$ is strictly stationary.
Furthermore, it is not $\alpha$-mixing in the sense of Rosenblatt \cite{bradley1985basic}, but it is $\tilde{\phi}$-dependent in the sense of Dedecker-Prieur~\cite{dedecker_prieur} (see also Caron-Dede~\cite{caron1}, Section $4.4$). Indeed, one can prove that the coefficients $\tilde{\phi}(k)$ of the chain $(Z_{i})_{i \geq 1}$ decrease geometrically \cite{dedecker_prieur}: $\tilde{\phi}(k) \leq 2^{-k}$.
Let now $Q_{0,\sigma^{2}}$ be the inverse of the cumulative distribution function of the law $\mathcal{N}(0, \sigma^{2})$. Let then:
\[\epsilon_{i} = Q_{0,\sigma^{2}}(Z_{i}).\]

The sequence $(\epsilon_{i})_{i \geq 1}$ is also a stationary Markov chain (as an invertible function of a stationary Markov chain), and one can easily check that its $\tilde{\phi}(k)$ coefficients are exactly equal to those of the sequence $(Z_{i})_{i \geq 1}$ (hence, $(\epsilon_{i})$ satisfies Hannan's condition~\eqref{C1}, see Section $5.1$ in~\cite{caron1}). By construction, $\epsilon_{i}$ is $\mathcal{N}(0, \sigma^{2})$-distributed, but the sequence $(\epsilon_{i})_{i \geq 1}$ is not a Gaussian process (otherwise it would be mixing in the sense of Rosenblatt).
Consequently Hannan's conditions are satisfied and the tests can be corrected as indicated above.
For the simulations, let us notice that the variance $\sigma^{2}$ is chosen equal to $25$.\\

The first model simulated with this error process is the following linear regression model, for all $i$ in $\{1, \ldots, n\}$:
\begin{equation}
Y_{i} = \beta_{0} + \beta_{1}(i^{2} + X_{i}) + \epsilon_{i},
\label{Model1}
\end{equation}
with $(X_{i})_{i \geq 1}$ a gaussian $AR(1)$ process (the variance is equal to 9), independent of the Markov chain $(\epsilon_{i})_{i \geq 1}$.
The coefficient $\beta_{0}$ is chosen equal to $3$. 

We test the hypothesis $H_{0}$: $\beta_{1} = 0$, against the hypothesis $H_{1}$: $\beta_{1} \neq 0$, thanks to the statistic $T_{j,n}$ defined above~\eqref{Student}.
The estimated level of the test will be studied for different choices of $n$ and $h_{n}$, which is linked to the number of covariance terms considered.
Under the hypothesis $H_{0}$, the same test is carried out $2000$ times. Then we look at the frequency of rejection of the test when we are under $H_{0}$, that is to say the estimated level of the test. Let us specify that we want an estimated level close to $5\%$.\\

$\bullet$ Case $\beta_{1} = 0$ and $h_{n} = 1$ (no correction): \\

\begin{center}
\begin{tabular}{|c|c|c|c|c|c|}
\hline
$n$ & 200 & 400 & 600 & 800 & 1000 \\
\hline
Estimated level & 0.203 & 0.195 & 0.183 & 0.205 & 0.202 \\
\hline
\end{tabular} \\
\end{center}

Here, since $h_{n} = 1$, we do not estimate any of the covariance terms. The result is that the estimated levels are too large. This means that the test will reject the null hypothesis too often. \\

\begin{figure}[h]
\begin{center}
\includegraphics[scale=0.42]{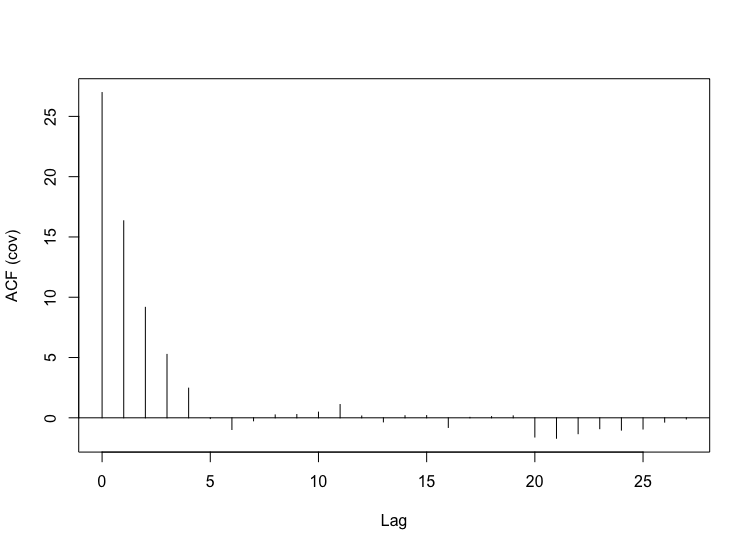} 
\end{center}
\caption{Empirical autocovariance of the residuals of Model~\eqref{Model1}.}
\label{fig1}
\end{figure}

The parameter $h_{n}$ may be chosen by analyzing the graph of the empirical autocovariances, Figure~\ref{fig1}. For this example, the shape of the empirical autocovariance suggests to keep only $4$ terms. This leads to choose $h_{n} =  5$. \\

$\bullet$ Case $\beta_{1} = 0$, $h_{n} = 5$: \\

\begin{center}
\begin{tabular}{|c|c|c|c|c|c|}
\hline
$n$ & 200 & 400 & 600 & 800 & 1000 \\
\hline
Estimated level & 0.0845 & 0.065 & 0.0595 & 0.054 & 0.053 \\
\hline
\end{tabular} \\
\end{center}

As suggested by the graph of the empirical autocovariances, the choice $h_{n} = 5$ gives better estimated levels than $h_{n} = 1$.
If one increases the size of the samples $n$, we are getting closer to the estimated level $5$\%. If $n = 2000$, the estimated level is around $0.05$.

Let us notice that even for $n$ moderately large ($n$ approximately $200$), it is much better to correct the test than not to do it. The estimated level goes from $20\%$ to $8.5\%$.\\

$\bullet$ Case $\beta_{1} = 0.00001$, $h_{n} = 5$: \\

In this example, $H_{0}$ is not satisfied. We choose $\beta_{1}$ equal to $0.00001$, and  we perform the same tests as above ($N=2000$) to estimate the power of the test.

\begin{center}
\begin{tabular}{|c|c|c|c|c|c|}
\hline
$n$ & 200 & 400 & 600 & 800 & 1000 \\
\hline
Estimated power & 0.1025 & 0.301 & 0.887 & 1 & 1 \\
\hline
\end{tabular} \\
\end{center}

As one can see, the estimated power is always greater than $0.05$, as expected.
Still as expected, the estimated power increases with the size of the samples. For $n = 200$, the power of the test is around $0.10$, and for $n = 800$, the power is around $1$.
As soon as $n = 800$, the test always rejects the $H_{0}$-hypothesis.\\

The second model considered is the following linear regression model, for all $i$ in $\{1, \ldots, n\}$:
\begin{equation}
Y_{i} = \beta_{0} + \beta_{1}(\log(i) + \sin(i) + X_{i}) + \beta_{2} i + \epsilon_{i}
\label{Model2}
\end{equation}

Here, we test the hypothesis $H_{0}$: $\beta_{1} = \beta_{2} = 0$ against $H_{1}$: $\beta_{1} \neq 0$ or $\beta_{2} \neq 0$, thanks to the statistic~$\Xi$~\eqref{Fisher}. The coefficient $\beta_{0}$ is equal to $3$, and we use the same simulation scheme as above. \\

$\bullet$ Case $\beta_{1} = \beta_{2} = 0$ and $h_{n} = 1$ (no correction): \\

\begin{center}
\begin{tabular}{|c|c|c|c|c|c|}
\hline
$n$ & 200 & 400 & 600 & 800 & 1000 \\
\hline
Estimated level & 0.348 & 0.334 & 0.324 & 0.3295 & 0.3285 \\
\hline
\end{tabular} \\
\end{center}

As for the first simulation, if $h_{n} = 1$ the test will reject the null hypothesis too often.\\

\begin{figure}[h]
\begin{center}
\includegraphics[scale=0.42]{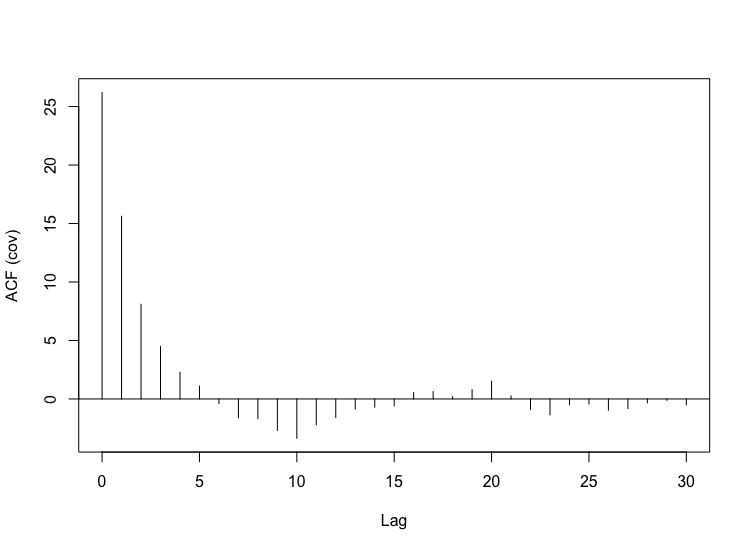} 
\end{center}
\caption{Empirical autocovariance of the residuals of Model~\eqref{Model2}.}
\label{fig2}
\end{figure}

As suggested by the graph of the estimated autocovariances Figure~\ref{fig2}, it suggests to keep only $5$ terms of covariances. Given the kernel~\eqref{kernel}, if we want to keep $5$ terms of covariances, we must choose a bandwidth equal to $h_{n} =  6.25$ (because $\frac{5}{0.8} = 6.25$). \\

$\bullet$ Case $\beta_{1} = \beta_{2} = 0$, $h_{n} = 6.25$: \\

\begin{center}
\begin{tabular}{|c|c|c|c|c|c|}
\hline
$n$ & 200 & 400 & 600 & 800 & 1000 \\
\hline
Estimated level & 0.09 & 0.078 & 0.066 & 0.0625 & 0.0595 \\
\hline
\end{tabular} \\
\end{center}

Here, we see that the choice $h_{n} = 6.25$ works well. For $n=1000$, the estimated level is around $0.06$. If $n = 2000$ and $h_{n} = 6.25$, the estimated level is around $0.05$.\\ 

$\bullet$ Case $\beta_{1} = 0.2$, $\beta_{2} = 0$, $h_{n} = 6.25$: \\

Now, we study the estimated power of the test. The coefficient $\beta_{1}$ is chosen equal to $0.2$ and $\beta_{2}$ is equal to $0$.

\begin{center}
\begin{tabular}{|c|c|c|c|c|c|}
\hline
$n$ & 200 & 400 & 600 & 800 & 1000 \\
\hline
Estimated power & 0.33 & 0.5 & 0.6515 & 0.776 & 0.884 \\
\hline
\end{tabular} \\
\end{center}

As expected, the estimated power increases with the size of the samples, and it is around $0.9$ when $n = 1000$.\\

\section{Proofs}

\subsection{Theorem~\ref{main_th}}

\begin{proof}

In this proof, we use the notations introduced in Section $2$ and Section $3$. We denote by $V(X)$ the matrix equal to $\mathbb{E} \left( D(n) (\hat{\beta} -\beta) (\hat{\beta} -\beta)^{t} D(n)^{t} \Big| X \right)$ and by $v_{j,l}$ its coefficients.\\

By the triangle inequality, we have for all $j, l$ in $\{1, \ldots, p\}$:
\[\left| c_{n,(j,l)} - c_{j,l} \right| \leq \left| v_{j,l} - c_{j,l} \right| + \left| c_{n,(j,l)} - v_{j,l} \right|.\]
Thanks to Hannan's Theorem~\ref{Hannan_th}, we already know that:
\begin{equation}
\lim_{n \rightarrow \infty} \mathbb{E} \left( \left| v_{j,l} - c_{j,l} \right| \Big| X \right) = 0, \quad a.s.
\label{main_th_eq1}
\end{equation}
Then it remains to prove that:
\begin{equation}
\lim_{n \rightarrow \infty} \mathbb{E} \left( \left| c_{n,(j,l)} - v_{j,l} \right| \Big| X \right) = 0, \quad a.s.
\label{main_th_eq2}
\end{equation}

The matrix $V(X)$ is equal to:
\[D(n) (X^{t}X)^{-1} X^{t} \Gamma_{n} X (X^{t}X)^{-1} D(n),\]
with $\Gamma_{n}$ defined in~\eqref{Grand_Gamma}, and the estimator $C_{n}$:
\[D(n) (X^{t}X)^{-1} X^{t} \widehat{\Gamma}_{n, h_{n}}^{\ast} X (X^{t}X)^{-1} D(n),\]
with $\widehat{\Gamma}_{n, h_{n}}^{\ast}$ defined in~\eqref{Gamma_tapered_star}.
Thanks to the convergence of $D_{n} (X^{t} X)^{-1} D_{n}$ to $R(0)^{-1}$, it is sufficient to consider the matrices: 
\[V' = D_{n}^{-1} X^{t} \Gamma_{n} X D_n^{-1},\]
and:
\[C_{n}' = D_{n}^{-1} X^{t} \widehat{\Gamma}_{n, h_{n}}^{\ast} X D_n^{-1}.\]
We know that $\Gamma_{n} = \sum_{k=-n+1}^{n-1} \gamma(k) J_{n}^{(k)}$ (see Remark~\ref{Rem2.2} for the definition of $J_{n}^{(k)}$). Thus we have for $V'$ and $C_{n}'$ the following decomposition:
\[D(n)^{-1} X^{t} \Gamma_{n} X D(n)^{-1} = \sum_{k=-n+1}^{n-1} \gamma(k) B_{k,n}\]
and:
\[D(n)^{-1} X^{t} \widehat{\Gamma}_{n, h_{n}}^{\ast} X D(n)^{-1} = \sum_{k=-n+1}^{n-1} K \left( \frac{k}{h_{n}} \right) \hat{\gamma}^{\ast}_{k} B_{k,n},\]
with:
\[B_{0,n} = D(n)^{-1} X^{t} X D(n)^{-1}\]
\[B_{k,n} = D(n)^{-1} X^{t} J_{n}^{(k)} X D(n)^{-1},\]
and:
\[\hat{\gamma}^{\ast}_{k} = \frac{1}{n} \sum_{j=1}^{n-|k|} \hat{\epsilon}_{j} \hat{\epsilon}_{j+|k|}. \]

Let us compute:
\[\left| c_{n,(j,l)}' - v_{j,l}' \right| = \left| \sum_{k=-n+1}^{n-1} \left( K \left( \frac{k}{h_{n}} \right) \hat{\gamma}_{k}^{\ast} - \gamma(k) \right) b_{j,l}^{k,n} \right|,\]
where $b_{j,l}^{k,n}$ is the coefficient $(j,l)$ of the matrix $B_{k,n}$. We shall show that:
\[\lim_{n \rightarrow \infty} \mathbb{E} \left( \left| \sum_{k=-n+1}^{n-1} \left( K \left( \frac{k}{h_{n}} \right) \hat{\gamma}_{k}^{\ast} - \gamma(k) \right) b_{j,l}^{k,n} \right| \Big| X \right) = 0, \quad a.s.\]

We recall that: 
\[f(\lambda) = \frac{1}{2 \pi} \sum_{k=-\infty}^{\infty} \gamma(k) e^{i k \lambda}, \qquad \gamma(k) = \int_{-\pi}^{\pi} e^{i k \lambda} f(\lambda) d \lambda,\]
where the coefficients $\gamma(k)$ are the Fourier coefficients of the spectral density $f(\lambda)$.
We have:
\[f_{n}^{\ast}(\lambda) = \frac{1}{2 \pi} \sum_{k=-n+1}^{n-1} K \left( \frac{k}{h_{n}} \right) \hat{\gamma}_{k}^{\ast} e^{i k \lambda}, \qquad K \left( \frac{k}{h_{n}} \right) \hat{\gamma}^{\ast}_{k} = \int_{-\pi}^{\pi} e^{i k \lambda} f_{n}^{\ast}(\lambda) d \lambda\]
and the coefficients $K \left( \frac{k}{h_{n}} \right) \hat{\gamma}^{\ast}_{k}$ are the Fourier coefficients of the spectral density's estimator $f_{n}^{\ast}(\lambda)$.
Let us define:
\[g_{n}(\lambda) = \frac{1}{2 \pi} \sum_{k=-n+1}^{n-1} e^{i k x} B_{k,n},\]
in such a way that the matrices $B_{k,n}$ are the Fourier coefficients of the function $g_{n}(\lambda)$:
\[B_{k,n} = \int_{-\pi}^{\pi} e^{i k \lambda} g_{n}(\lambda) d \lambda.\]
Consequently we can deduce that:
\[\sum_{k=-n+1}^{n-1} \left( K \left( \frac{k}{h_{n}} \right) \hat{\gamma}^{\ast}_{k} - \gamma(k) \right) B_{k,n} = \int_{-\pi}^{\pi} \left( f_{n}^{\ast}(\lambda) - f(\lambda) \right) g_{n}(\lambda) (d \lambda).\]
Thus it remains to prove that, for all $j, l$ in $\{1, \ldots, p\}$:
\[\lim_{n \rightarrow \infty} \mathbb{E} \left( \left| \int_{-\pi}^{\pi} \left( f_{n}^{\ast}(\lambda) - f(\lambda) \right) [g_{n}(\lambda)]_{j,l} d \lambda \right| \Big| X \right) = 0, \quad a.s.\]

We have:
\begin{eqnarray*}
\mathbb{E} \left( \left| \int_{-\pi}^{\pi} \left( f_{n}^{\ast}(\lambda) - f(\lambda) \right) [g_{n}(\lambda)]_{j,l} d \lambda \right| \Big| X \right) 
&\leq & \mathbb{E} \left( \int_{-\pi}^{\pi} \left| f_{n}^{\ast}(\lambda) - f(\lambda) \right| \left| [g_{n}(\lambda)]_{j,l} \right| d \lambda \Big| X \right) \\
&\leq & \int_{-\pi}^{\pi} \left| [g_{n}(\lambda)]_{j,l} \right|  \mathbb{E} \left( \left| f_{n}^{\ast}(\lambda) - f(\lambda) \right| \Big| X \right) d \lambda,
\end{eqnarray*}
because $\left[ g_{n}(\lambda) \right]_{j,l}$ is measurable with respect to the $\sigma$-algebra generated by the design $X$. Then:
\begin{multline*}
\int_{-\pi}^{\pi} \left | [g_{n}(\lambda)]_{j,l} \right| \mathbb{E} \left( \left| f_{n}^{\ast}(\lambda) - f(\lambda) \right|  \Big| X \right) d \lambda \\ 
\leq \sup_{\lambda \in [-\pi,\pi]} \mathbb{E} \left( \left| f_{n}^{\ast}(\lambda) - f(\lambda) \right| \Big| X \right) \int_{-\pi}^{\pi} \left | [g_{n}(\lambda)]_{j,l} \right| d \lambda.
\end{multline*}
Theorem $3.1$ of our paper~\cite{caron1} states that:
\[\lim_{n \rightarrow \infty} \sup_{\lambda \in [-\pi,\pi]} \left\| f_{n}^{\ast}(\lambda) - f(\lambda) \right \|_{\mathbb{L}^{1}} = 0,\]
for a fixed design $X$ and for the particular kernel defined by: $K(x) = \mathds{1}_{\{|x| \leq 1\}} + (2 - |x|) \mathds{1}_{\{1 \leq |x| \leq 2\}}$.
But a quick look to the proof of this theorem suffices to see that this result is available for any design $X$, conditionally to $X$: \[\lim_{n \rightarrow \infty} \sup_{\lambda \in [-\pi,\pi]} \mathbb{E} \left( \left| f_{n}^{\ast}(\lambda) - f(\lambda) \right| \Big| X \right) = 0, \quad a.s.\]
Furthermore, this result is still available for all kernel $K$ verifying the conditions at the beginning of Section $3$.

Thus it remains to find a bound for:
\[\int_{-\pi}^{\pi} \left | [g_{n}(\lambda)]_{j,l} \right| d \lambda.\]
Let us recall (see~\eqref{Bkn}) that the matrices $B_{k,n}$ are equal to, for all $k$ in $\{-n+1, \dots, n-1\}$:
\[B_{k,n} = [\hat{\rho}_{j,l}(k)], \quad \text{if} \ k \geq 0, \qquad B_{k,n} = [\hat{\rho}_{j,l}(-k)], \quad \text{if} \ k \leq -1.\]

By definition we have:
\begin{equation}
\hat{\rho}_{j,l}(k) = \frac{\hat{\gamma}_{j,l}(k)}{\sqrt{\hat{\gamma}_{j,j}(0) \hat{\gamma}_{l,l}(0)}}.
\label{rho_chapeau}
\end{equation}
For a multivariate time series, let us recall that the cross-periodogram is defined by, for all $j$, $l$ in $\{1, \ldots, p\}$~\cite{brockwell2013time}:
\begin{equation}
[I_{n}(\lambda)]_{j,l} = \frac{1}{2 \pi} \sum_{k=-n+1}^{n-1} \hat{\gamma}_{j,l}(k) e^{i k \lambda}.
\label{cross-per}
\end{equation}
Combining~\eqref{rho_chapeau} and~\eqref{cross-per}, the function $g_{n}(\lambda)$ is equal to, for all $j$, $l$ in $\{1, \ldots, p\}$:
\[[g_{n}(\lambda)]_{j,l} = \frac{[I_{n}(\lambda)]_{j,l}}{\sqrt{\hat{\gamma}_{j,j}(0) \hat{\gamma}_{l,l}(0)}} = \frac{1}{2 \pi} \sum_{k=-n+1}^{n-1} b_{j,l}^{k,n} e^{i k \lambda}.\]
Then using the definition of the coherence~\cite{brockwell2013time}, we get:
\begin{multline*}
\left| [g_{n}(\lambda)]_{j,l} \right| = \frac{\left| [I_{n}(\lambda)]_{j,l} \right|}{\sqrt{\hat{\gamma}_{j,j}(0) \hat{\gamma}_{l,l}(0)}} \leq \sqrt{\frac{[I_{n}(\lambda)]_{j,j} [I_{n}(\lambda)]_{l,l}}{\hat{\gamma}_{j,j}(0) \hat{\gamma}_{l,l}(0)}} \\ 
\leq \sqrt{[g_{n}(\lambda)]_{j,j} [g_{n}(\lambda)]_{l,l}} \leq \frac{1}{2} [g_{n}(\lambda)]_{j,j} + \frac{1}{2} [g_{n}(\lambda)]_{l,l}.
\end{multline*}
Consequently, we have:
\[\int_{-\pi}^{\pi} \left | [g_{n}(\lambda)]_{j,l} \right| d \lambda \leq \frac{1}{2} \int_{-\pi}^{\pi} [g_{n}(\lambda)]_{j,j} d \lambda + \frac{1}{2} \int_{-\pi}^{\pi} [g_{n}(\lambda)]_{l,l} d \lambda \leq \frac{1}{2} [B_{0,n}]_{j,j} + \frac{1}{2} [B_{0,n}]_{l,l} \leq 1,\]
because $[B_{0,n}]_{j,j} = \hat{\rho}_{j,j}(0) = 1$ and $[B_{0,n}]_{l,l} = \hat{\rho}_{l,l}(0) = 1$.

We deduce that, for all $j$, $l$ in $\{1, \ldots, p\}$:
\begin{multline*}
\mathbb{E} \left( \left| \int_{-\pi}^{\pi} \left( f_{n}^{\ast}(\lambda) - f(\lambda) \right) [g_{n}(\lambda)]_{j,l} d \lambda \right| \Big| X \right) \\
\leq \sup_{\lambda \in [-\pi,\pi]} \mathbb{E} \left( \left| f_{n}^{\ast}(\lambda) - f(\lambda) \right| \Big| X \right) \int_{-\pi}^{\pi} \left | [g_{n}(\lambda)]_{j,l} \right| d \lambda \\
\leq \sup_{\lambda \in [-\pi,\pi]} \mathbb{E} \left( \left| f_{n}^{\ast}(\lambda) - f(\lambda) \right| \Big| X \right).
\end{multline*}
Since we know that:
\[\lim_{n \rightarrow \infty} \sup_{\lambda \in [-\pi,\pi]} \mathbb{E} \left( \left| f_{n}^{\ast}(\lambda) - f(\lambda) \right| \Big| X \right) = 0, \quad a.s.\]
we have proved that, for all $j$, $l$ in $\{1, \ldots, p\}$:
\[\lim_{n \rightarrow \infty} \mathbb{E} \left( \left| \int_{-\pi}^{\pi} \left( f_{n}^{\ast}(\lambda) - f(\lambda) \right) [g_{n}(\lambda)]_{j,l} d \lambda \right| \Big| X \right) = 0, \quad a.s.\]

\end{proof}

\subsection{Corollary~\ref{cor_main_th}}

\begin{proof}

We want to prove that, for all $j$, $l$ in $\{1, \ldots, p\}$, $c_{n,(j,l)}$ converges in probability to $c_{j,l}$ as $n$ tends to infinity, that is, for all $\epsilon > 0$:
\[\mathbb{E} \left( \mathds{1}_{| c_{n,(j,l)} - c_{j,l} | > \epsilon} \right) \xrightarrow[n \rightarrow \infty]{} 0.\]
We have:
\[\mathbb{E} \left( \mathds{1}_{| c_{n,(j,l)} - c_{j,l} | > \epsilon} \right) = \mathbb{E} \left( \mathbb{E} \left( \mathds{1}_{| c_{n,(j,l)} - c_{j,l} | > \epsilon} \Big| X \right) \right).\]
Thanks to Theorem~\ref{main_th} and to Markov's inequality, we have almost surely:
\[\mathbb{E} \left( \mathds{1}_{| c_{n,(j,l)} - c_{j,l} | > \epsilon} | X \right) \leq \frac{\mathbb{E} \left( \left| c_{n,(j,l)} - c_{j,l} \right| \Big| X \right)}{\epsilon} \xrightarrow[n \rightarrow \infty]{} 0.\]
Then, using the dominated convergence theorem, we get:
\[\mathbb{E} \left( \mathbb{E} \left( \mathds{1}_{| c_{n,(j,l)} - c_{j,l} | > \epsilon} \Big| X \right) \right) \xrightarrow[n \rightarrow \infty]{} 0.\]

\end{proof}

\newpage

\bibliographystyle{abbrv}
\bibliography{Caron_bib}
\addcontentsline{toc}{chapter}{Bibliographie}

\end{document}